\documentclass[11pt]{article}

\usepackage{latexsym,amsmath,amssymb}

\begin{document}

\centerline{\large AN IDENTITY OF ANDREWS AND A NEW METHOD FOR
THE}

\centerline{\large RIORDAN ARRAY PROOF OF COMBINATORIAL
IDENTITIES}

\vspace{1.5cm}

\centerline{\large Eduardo H.\ M.\ Brietzke}

\centerline{\small Instituto de Matem\'atica -- UFRGS}

\centerline{\small Caixa Postal 15080}

\centerline{\small 91509--900 Porto Alegre, RS, Brazil}

\centerline{\small email: brietzke@mat.ufrgs.br}

\vspace{0.7cm}

\begin{abstract}
We consider an identity relating Fibonacci numbers to Pascal's
triangle discovered by \mbox{G.\ E.\ Andrews}. Several authors
provided proofs of this identity, most of them rather involved or
else relying on sophisticated number theoretical arguments. We
present a new proof, quite simple and based on a Riordan array
argument. The main point of the proof is the construction of a new
Riordan array from a given Riordan array, by the elimination of
elements. We extend the method and as an application we
obtain other identities, some of which are new. An important
feature of our construction is that it establishes a nice
connection between the generating function of the $A-$sequence of
a certain class of Riordan arrays and hypergeometric functions.
\end{abstract}


\vspace{0.5cm}

\section{Introduction}

In this article we provide a new proof of an identity of Andrews,
based on Riordan arrays. Several authors have already proved this
identity using different types of argument (references are given
below). Our reason to give a further proof is that we believe our
idea is new and interesting on its own.

In our approach we establish a nice connection between the
generating function of the $A-$sequence of a certain class of
Riordan arrays and hypergeometric functions. This new connection
with hypergeometric functions is probably in itself interesting
and is one of the main features of this work. Our method involves
constructing a new \mbox{Riordan} array from a given Riordan array
by eliminating entire rows and parts of the remaining rows. In the
proof of the identity of Andrews, this construction is applied to
Pascal's triangle, but for the sake of illustrating the usefulness
of our method, we make additional applications to Pascal's
triangle as well as to other Riordan arrays, for example Catalan's
triangle, obtaining a few more identities.

As a generalization of Pascal's, Catalan's, Motzkin's, and other
triangles, D.\ G.\ Rogers introduced in 1978 (\cite{20}) the
concept of renewal array, which was further generalized to Riordan
array by Shapiro et al. in 1991 (\cite{22}). Among other
applications Riordan arrays turned out to be an extremely powerful
tool in dealing with combinatorial identities. R.\ Sprugnoli in
\cite{15} used Riordan arrays to find several combinatorial sums
in closed form and also to determine their asymptotic value. For
additional applications of Riordan arrays to the evaluation in
closed form of sums involving binomial, Stirling, Bernoulli, and
harmonic numbers, see \cite{28}. The Riordan array technique has
also been employed to show that two combinatorial sums are
equivalent, regardless of whether they have a closed form
expression or not (see \cite{23}). An important problem that has
occupied mathematicians for a long time is the inversion of
combinatorial sums (see \cite{25}). The concept of Riordan array
provided a powerful tool to prove a large class of inversions
(see, for example, \cite{26} and \cite{24}).

The paper is organized as follows. In the introduction we recall
some basic results needed in the sequel. In Section 2 we develop
our method of extracting new Riordan arrays from a given one. We
also establish a connection between a certain class of Riordan
arrays and hypergeometric functions. The ideas of Section 2 are
applied in Section 3 to give a new proof of the identity of
Andrews. Section 4 is devoted to additional applications of our
method. As an illustration of our ideas, further identities are
obtained. Identities of this type can often be proved directly,
using generating functions and Lagrange's Inversion Formula. To
show how this can be done, in Section 5 we give a direct proof of
one of the identities obtained previously.

We begin by recalling Lagrange's Inversion Theorem, which is an
important element needed in our study. Several forms of Lagrange's
Inversion Formula exist (see \cite{30}). We summarize some of them
below.

\vspace{0.6cm}

\noindent {\bf Theorem 1.1} (Lagrange's Inversion Theorem
\cite{30}) {\em Suppose that a formal power series $w=w(t)$ is
implicitly defined by the relation $w=t\phi(w)$, where $\phi(t)$
is a formal power series such that $\phi(0)\neq0$. Then,
\begin{equation}
[t^n]\big(w(t)\big)^k=\frac{k}{n}\,[t^{n-k}]\big(\phi(t)\big)^n.
\label{e:lag1}
\end{equation}
Equivalently, for any formal power series $F(t)$,
\begin{equation}
[t^n]F(w(t))=\frac{1}{n}[t^{n-1}]F'(t)\big(\phi(t)\big)^n.
\label{e:lag3}
\end{equation}
In terms of generating functions,}
\begin{equation}
\mathcal{G}\big([t^n]F(t)(\phi(t))^n\big)=
\bigg[\frac{F(w)}{1-t\phi'(w)}\,\bigg|\,w=t\phi(w)\bigg].
\label{e:lag5}
\end{equation}

\vspace{0.4cm}

The above notation, i.e., $[f(w)|w=g(t)]$, means replacing $w$ by
$g(t)$ in $f(w)$, and given any sequence $(b_n)$,
$\mathcal{G}(b_n)$ stands for its generating function
$\mathcal{G}(b_n)=\sum_{n=0}^\infty b_nt^n$.

A Riordan array is an infinite lower triangular array
$D=\{d_{n,k}\}_{n,k\ge0}$ defined by a pair of formal power series
$D=\big(d(t),h(t)\big)$, for which
\begin{equation}
d_{n,k}=[t^n]\,d(t)\big(th(t)\big)^k, \quad \forall n,k\ge0.
\label{e:1}
\end{equation}
Here $[t^n]g(t)$ denotes the coefficient of $t^n$ in $g(t)$.
Pascal's triangle is an example of a Riordan array. In this case,
$d(t)=h(t)=1/(1-t)$ and $d_{n,k}= \binom{n}{k}$. One of the main
results of the theory of Riordan arrays is the following theorem.

\vspace{0.6cm}

\noindent {\bf Theorem 1.2} (\cite{15},Theorem 1.1) {\em Let
$D=\big(d(t),h(t)\big)$ be a Riordan array and $f(t)=\sum f_kt^k$
a formal power series. Then,}
\begin{equation}
\sum_{k=0}^\infty f_kd_{n,k}=[t^n]\,d(t)f\big(th(t)\big).
\label{e:2}
\end{equation}

\vspace{0.4cm}

A Riordan array $D=\big(d(t),h(t)\big)$ is called proper if
$h_0=h(0)\neq0$. In \cite{20}, D.\ Rogers pointed out that proper
Riordan arrays can be alternately characterized by a pair
$d(t)=\sum_nd_{n,0}t^n$, the generating function of the first
column, and $A(t)=\sum a_kt^k$, the generating function of the
$A-$sequence, such that
\begin{equation}
d_{n+1,k+1}=a_0d_{n,k}+a_1d_{n,k+1}+a_2d_{n,k+2}+\cdots, \quad
\forall n,k\ge0.
\label{e:3}
\end{equation}
If $D=\big(d(t),h(t)\big)$ is a proper Riordan array, then
ord$\big((th(t))^k\big)=k$, for every $k$, where for any non-zero
formal power series $g(t)$ the order ord$\big(g(t)\big)$ is the
index of the first non-zero coefficient of $g(t)$. Therefore,
there exists a unique sequence $(a_k)$, called the $A-$sequence of
the Riordan array, such that
\[
h(t)=a_0+a_1th(t)+a_2\big(th(t)\big)^2+\cdots
\]
[i.e., the formal power series $A(t)=a_0+a_1t+a_2t^2+\cdots$,
refered to as the generating function of the $A-$sequence, is such
that $h(t)=A\big(th(t)\big)$]. Multiplying by
$d(t)\big(th(t)\big)^k$, we obtain
\[
t^{-1}d(t)(th(t))^{k+1}=
a_0d(t)(th(t))^k+a_1d(t)(th(t))^{k+1}+a_2d(t)(th(t))^{k+2}+\cdots
\]
and applying $[t^n]$ to both sides, \eqref{e:3} follows.

The converse is also true and we state it as the following
theorem.

\vspace{0.6cm}

\noindent {\bf Theorem 1.3} (\cite{15}, Theorem 1.3) {\em Let
$D=\{d_{n,k}\}_{n\ge k\ge0}$ be an infinite triangle such that
$d_{0,0}\neq0$ and for which \eqref{e:3} holds for some sequence
$(a_k)$ with $a_0\neq0$. Then $D$ is a proper Riordan array
$\big(d(t),h(t)\big)$, where $d(t)=\sum_{n=0}^\infty d_{n,0}t^n$
is the generating function of the first column and $h(t)$ is the
unique solution of
\begin{equation}
h(t)=A\big(th(t)\big),
\label{e:4}
\end{equation}
for $A(t)=\mathcal{G}(a_k)=\sum a_kt^k$ the generating funtion of
the sequence $(a_k)$. Moreover,}
\begin{equation}
[t^{n-1}]h(t)=\frac{1}{n}[t^{n-1}]\big(A(t)\big)^n.
\label{e:454545}
\end{equation}

\vspace{0.3cm}

\noindent {\bf Proof.} Since $a_0\neq0$, by Lagrange's Inversion
Formula \eqref{e:lag1} with $w:=th(t)$, $k=1$ and $\phi=A$,
\eqref{e:4} defines a unique formal power series $h(t)$, for which
\eqref{e:454545} holds. We have to verify that given $n$,
\[
d_{n,k}=[t^n]d(t)\big(th(t)\big)^k,
\]
for all $k$. By induction, suppose this holds for some $n$. Then,
\[
\begin{split}
d_{n+1,k+1}&=\sum_{j\ge0}a_jd_{n,k+j}=
\sum_{j\ge0}a_j[t^n]d(t)\big(th(t)\big)^{k+j}=
[t^n]d(t)\big(th(t)\big)^kA\big(th(t)\big)\\
&=[t^n]d(t)\big(th(t)\big)^kh(t)
=[t^{n+1}]d(t)\big(th(t)\big)^{k+1}.\rule{0.cm}{0.6cm}
\qquad\square
\end{split}
\]

\vspace{0.6cm}

\section{Construction of a new Riordan array}

\vspace{0.3cm}

We now describe a process of obtaining new Riordan arrays from a
given Riordan array, which corresponds to eliminating rows from
the original array, eliminating the first elements from the
remaining rows, and shifting them to the left. For a fixed $p$ we
keep one of every $p$ rows.

\vspace{0.5cm}

\noindent {\bf Theorem 2.1} {\em Given a proper Riordan array
$\{d_{n,k}\}_{n,k\ge0}$, for any integers $p\ge2$ and $r\ge0$,
$\tilde{d}_{n,k}=d_{pn+r,(p-1)n+r+k}$ $\,(n,k\ge0)$ defines a new
Riordan array. Moreover, the generating function of the
$A-$sequence of the new array is $\big(A(t)\big)^p$, where $A(t)$
is the generating function of the $A-$sequence of the given
Riordan array.}

\vspace{0.5cm}

\noindent {\bf Proof.} Let $A(t)=a_0+a_1t+a_2t^2+\cdots$ be as in
\eqref{e:3}. If $p=2$, for $r=0,1$ we have
\[
\begin{split}
\tilde{d}_{n+1,k+1}&=d_{2n+2+r,\,n+k+2+r}= \sum_{i=0}^\infty
a_id_{2n+1+r,\,n+k+i+1+r} \\
&=\sum_{i=0}^\infty\sum_{j=0}^\infty a_ia_jd_{2n+r,\,n+k+i+j+r}
\end{split}
\]
and, therefore,
\[
\tilde{d}_{n+1,k+1}=
\sum_{\nu=0}^\infty\sum_{i=0}^{\nu}a_ia_{r-i}d_{2n+r,n+k+\nu+r},
\]
i.e.,
\[
\tilde{d}_{n+1,k+1}=
\sum_{\nu=0}^\infty\sum_{i=0}^{\nu}a_ia_{r-i}\tilde{d}_{n,k+\nu}.
\hspace{0.5cm}
\]
Hence,
\[
\tilde{d}_{n+1,k+1}=\sum_{r=0}^\infty b_r\tilde{d}_{n,k+r}, \qquad
{\rm where} \quad \sum_{k=0}^\infty b_kt^k=\big(A(t)\big)^2.
\]
By Theorem \ref{t:2}, $\{\tilde{d}_{n,k}\}_{n,k\ge0}$ is a Riordan
array and $B(t)=\big(A(t)\big)^2$ is the generating function of
its $A-$sequence. If $p\ge3$ an iteration of the argument applies.
$\qquad\square$

\vspace{0.6cm}

For example, beginning with Pascal's triangle, for $p=3$ and
$r=1$, we obtain the Riordan array

\begin{center}

\begin{tabular}{ccccccc}
    1 &      &      &     &     &     &    \\
    4 &    1 &      &     &     &     &    \\
   21 &    7 &    1 &     &     &     &    \\
  120 &   45 &   10 &   1 &     &     &    \\
  715 &  286 &   78 &  13 &   1 &     &    \\
 2002 & 1001 &  560 & 120 &  16 &  1  &    \\
27132 &11628 & 3876 & 969 & 171 & 19  & \ 1 \\
\end{tabular}
\centerline{$\cdots \cdots$}

\end{center}

\noindent in which $d_{n,k}=\binom{3n+1}{2n+k+1}$ and, by Theorem
2.1,
\[
\hspace{2.cm}
d_{n+1,k+1}=d_{n,k}+3d_{n,k+1}+3d_{n,k+2}+d_{n,k+3}.
\]
In what follows we use the generalized hypergeometric series,
defined by
\[
\,_pF_q\bigg(
\begin{array}{c}
a_1,\ldots,a_p  \\[0pt]
c_1,\ldots,c_q
\end{array}
\bigg|\:t \bigg)=\sum_{n=0}^\infty\,
\frac{(a_1)_n\cdots(a_p)_n}{(c_1)_n\cdots(c_q)_n}\cdot\frac{t^n}{n!}\,,
\]
where $(a)_n$ stands for the Pochhammer symbol
\[
(a)_n=\frac{\Gamma(a+n)}{\Gamma(a)}=\left\{
\begin{array}{ll}
1 ,  & \mbox{if } n=0 \\
a(a+1)\cdots(a+n-1) , & \mbox{if } n\ge1
\end{array}
\right.
\]
The hypergeometric series is characterized by the fact that its
constant term is 1 and, setting
$A_n=\frac{(a_1)_n\cdots(a_p)_n}{(c_1)_n\cdots(c_q)_n}\,t^n$, the
ratio of consecutive terms is
\[
\frac{A_{n+1}}{A_n}=
\frac{(a_1+n)\cdots(a_p+n)}{(c_1+n)\cdots(c_q+n)}
\cdot\frac{t}{n+1}\,.
\]
In order to apply Theorem 2.1 we need the following result, which
establishes an interesting connection between Riordan arrays and
hypergeometric functions.

\vspace{0.6cm}

\noindent {\bf Theorem 2.2} {\em If the generating function of the
$A-$sequence of a proper Riordan array is $A(t)=(1+t)^q$, with
$q\in\mathbb{N}$, then $h$ is the hypergeometric function
\begin{equation}
h(t)=\,_{q}F_{q-1}\left(
\begin{array}{c}
\frac{q}{q},\frac{q+1}{q},\ldots,\frac{2q-1}{q}\\[4pt]
\frac{q+1}{q-1},\frac{q+2}{q-1},\ldots,\frac{2q-1}{q-1}
\end{array}
\Bigg| \: \frac{q^qt}{(q\!-\!1)^{q-1}}\right),
\label{e:23501}
\end{equation}
also given by
\begin{equation}
h(t)=\sum_{n=1}^\infty\frac{1}{(q-1)n+1}\binom{qn}{n}t^{n-1}=
\sum_{n=1}^\infty\frac{1}{qn+1}\binom{qn+1}{n}t^{n-1}.
\label{e:23421}
\end{equation}
Moreover,
\begin{equation}
\big(h(t)\big)^s=\,_{q}F_{q-1}\left(
\begin{array}{c}
\frac{sq}{q},\frac{sq+1}{q},\ldots,\frac{(s+1)q-1}{q}\\[4pt]
\frac{sq+1}{q-1},\frac{sq+2}{q-1},\ldots,\frac{(s+1)q-1}{q-1}
\end{array}
\Bigg| \frac{q^qt}{(q\!-\!1)^{q-1}}\right)
\label{e:23502}
\end{equation}
for every $s\in\mathbb{R}$. Consequently,
\[
\hspace{0.2cm} \big(h(t)\big)^s=\sum_{n=0}^\infty
\frac{qs}{(q-1)n+qs}\binom{q(n+s)-1}{n}t^n
\]
and, therefore,}
\begin{equation}
[t^j]\big(th(t)\big)^s=[t^{j-s}]\big(h(t)\big)^s=
\frac{qs}{(q-1)j+s}\binom{qj-1}{j-s}.
\label{e:1912}
\end{equation}

\vspace{0.5cm}

\noindent {\bf Proof.} If $A(t)=(1+t)^q$, by \eqref{e:454545},
Theorem 1.3,
\[
[t^{n-1}]h(t)=\frac{1}{n}[t^{n-1}](1+t)^{qn}=
\frac{1}{n}\binom{qn}{n-1}=
\frac{(qn)!}{\big((q\!-\!1)n+1\big)!\,n!}\,.
\]
Therefore,
\begin{equation}
h(t)=\sum_{n=1}^\infty\frac{1}{(q-1)n+1}\binom{qn}{n}t^{n-1}=
\sum_{n=0}^\infty\frac{1}{(q-1)n+q}\binom{qn+q}{n+1}t^n.
\label{e:15h}
\end{equation}
The series \eqref{e:15h} is hypergeometric as its constant
term is 1 and, setting
$A_n=\frac{(qn+q)!}{(n+1)!((q-1)n+q)!}t^n$, the ratio of
consecutive terms is
\[
\begin{split}
\frac{A_{n+1}}{A_n}&=\frac{(qn+q+1)(qn+q+2)\cdots(qn+2q)}
{\big((q\!-\!1)n\!+\!q\!+\!1\big)\big((q\!-\!1)n\!+\!q\!+\!2)\big)\cdots
\big((q\!-\!1)n\!+\!2q\!-\!1\big)}\!\cdot\!\frac{t}{n+2}
\\
\rule{0.cm}{1.cm}
&=\frac{\big(n+\frac{q+1}{q}\big)\big(n+\frac{q+2}{q}\big)\cdots\big(n+\frac{2q}{q}\big)}
{\big(n+\frac{q+1}{q-1}\big)\big(n+\frac{q+2}{q-1}\big)\cdots\big(n+\frac{2q-1}{q-1}\big)}
\!\cdot\!\frac{n+1}{n+2}\!\cdot\!\frac{t}{n+1}\!\cdot\!\frac{q^q}{(q-1)^{q-1}}
\\
\rule{0.cm}{1.cm}
&=\frac{\big(n+\frac{q}{q}\big)\big(n+\frac{q+1}{q}\big)\cdots\big(n+\frac{2q-1}{q}\big)}
{\big(n+\frac{q+1}{q-1}\big)\big(n+\frac{q+2}{q-1}\big)\cdots\big(n+\frac{2q-1}{q-1}\big)}
\!\cdot\!\frac{1}{n+1}\!\cdot\!\frac{q^qt}{(q-1)^{q-1}}\, .
\end{split}
\]
Thus, \eqref{e:23501} follows immediately.

\noindent Note that $h(t)=\frac{\mathcal{B}_q(t)-1}{t}$, where
$\mathcal{B}_q$ is the generalized binomial series, given by
\[
\hspace{0.3cm}
\mathcal{B}_q(t)=\sum_{n=0}^\infty\frac{1}{qn+1}\binom{qn+1}{n}t^n,
\]
for which we have
\begin{equation}
\big(\mathcal{B}_q(t)\big)^r=\sum_{n=0}^\infty\frac{r}{qn+r}\binom{qn+r}{n}t^n,
\hspace{0.3cm} \label{e:2352}
\end{equation}
for any real $r$ (see \cite{19}, page 201, and also \cite{30},
where in Theorem 2.1 this is obtained from Lagrange's Inversion
Theorem, using \eqref{e:lag5}). From \eqref{e:2352}, by the same
argument used for the series \eqref{e:15h}, it follows that, for
any $r\in\mathbb{R}$,
\begin{equation}
\big(\mathcal{B}_q(t)\big)^r= \,_{q}F_{q-1}\left(
\begin{array}{c}
\frac{r}{q},\frac{r+1}{q},\ldots,\frac{r+q-1}{q}\\[4pt]
\frac{r+1}{q-1},\frac{r+2}{q-1},\ldots,\frac{r+q-1}{q-1}
\end{array}
\Bigg| \: \frac{q^qt}{(q\!-\!1)^{q-1}}\right).
\label{e:23521}
\end{equation}
On the other hand, since
$h(t)=A\big(th(t)\big)=\big(1+th(t)\big)^q=\big(\mathcal{B}_q(t)\big)^q$,
replacing $r$ by $qs$ in \eqref{e:23521}, we have \eqref{e:23502}
for any $s\in\mathbb{R}$ and, therefore, \eqref{e:1912} holds.
$\qquad\square$

\vspace{0.5cm}

\noindent {\bf Remark.} From \eqref{e:23521}, we obtain the
remarkable identity

\[
\left[\,_{q}F_{q-1}\left(
\begin{array}{c}
\frac{1}{q},\frac{2}{q},\ldots,\frac{q}{q}\\[4pt]
\frac{2}{q-1},\frac{3}{q-1},\ldots,\frac{q}{q-1}
\end{array}
\Bigg|\:\frac{q^qt}{(q-1)^{q-1}}\right)\right]^r \hspace{2.5cm}
\]
\begin{equation}
\hspace{1.2cm}= \,_{q}F_{q-1}\left(
\begin{array}{c}
\frac{r}{q},\frac{r+1}{q},\ldots,\frac{r+q-1}{q}\\[4pt]
\frac{r+1}{q-1},\frac{r+2}{q-1},\ldots,\frac{r+q-1}{q-1}
\end{array}
\Bigg|\:\frac{q^qt}{(q\!-\!1)^{q-1}}\right)
\label{e:235213}
\end{equation}

\vspace{0.1cm}

\noindent for powers of a hypergeometric function, which is
essentially contained in (5.60) of \cite{19} but is not explicitly
stated in the literature. Indeed, the literature does not refer to
many instances in which a product of hypergeometric functions is
also hypergeometric. In Section 5 we obtain two more identities
involving a product of hypergeometric functions.

\vspace{0.6cm}

\section{The Identities of Andrews}

\vspace{0.2cm}

We now apply the procedure described in Theorem 2.1 of extracting
new Riordan arrays from a given one to provide a new proof of some
identities obtained by G.\ E.\ Andrews in \cite{1}, namely
\begin{equation}
F_n=\sum_{k=-\infty}^\infty(-1)^k \,
\binom{n-1}{\big\lfloor\frac{1}{2}(n-1-5k)\big\rfloor}
\hspace{0.15cm} \label{e:a1}
\end{equation}
and
\begin{equation}
F_n=\sum_{k=-\infty}^\infty(-1)^k \, \binom{n}
{\big\lfloor\frac{1}{2}(n-1-5k)\big\rfloor}, \label{e:a2}
\end{equation}
where $(F_n)$ is the sequence of Fibonacci numbers, defined by
$F_0=0$, $F_1=1$, and $F_{n+2}=F_{n+1}+F_{n}$. Different proofs of
\eqref{e:a1} and \eqref{e:a2} were given by H.\ Gupta in \cite{3}
and by M.\ D.\ Hirschhorn in \cite{4} and \cite{5}. They are all
rather involved, though elementary, and they are specifically
designed to deal with the case of Pascal's triangle. As indicated
in \cite{3} and \cite{4}, identities \eqref{e:a1} and \eqref{e:a2}
are equivalent to
\begin{equation}
F_{2n+1}=\sum_{j=-\infty}^\infty
\biggl[\binom{2n+1}{n-5j}-\binom{2n+1}{n-5j-1}\biggr],
\label{e:a3}
\end{equation}
\begin{equation}
F_{2n+2}=\sum_{j=-\infty}^\infty
\biggl[\binom{2n+2}{n-5j}-\binom{2n+2}{n-5j-1}\biggr]
\hspace{0.2cm}
\label{e:121}
\end{equation}
and
\begin{equation}
F_{2n+2}=\sum_{j=-\infty}^\infty
\biggl[\binom{2n+1}{n-5j}-\binom{2n+1}{n-5j-2}\biggr],
\label{e:a5}
\end{equation}
\begin{equation}
F_{2n+1}=\sum_{j=-\infty}^\infty
\biggl[\binom{2n}{n-5j}-\binom{2n}{n-5j-2}\biggr],
\label{e:a6}
\end{equation}
respectively. In \cite{2} G.\ E.\ Andrews proves these identities
in the context of identities of the Rogers--Ramanujan type (see
also \cite{6}). In \cite{13}, identities \eqref{e:a3} through
\eqref{e:a6}, as well as several other similar identities for
trinomial coefficients and Catalan's triangle, have been proved in
a very elementary and direct way.

We now prove \eqref{e:121} to illustrate how identities
\eqref{e:a3} through \eqref{e:a6} can be obtained by a Riordan
array technique. Replacing $n$ by $n-1$ in \eqref{e:121}, it
suffices to show that
\begin{equation}
F_{2n}=\sum_{j=-\infty}^\infty
\biggl[\binom{2n}{n-5j-1}-\binom{2n}{n-5j-2}\biggr].
\hspace{0.2cm}
\label{e:122}
\end{equation}
We start with a visualization of Pascal's triangle in which
alternate rows have been removed and only non-vanishing binomial
numbers are represented:

\begin{center}

\begin{tabular}{ccccccccccccc}

       &    &    &     &            &           &   1 &     &
&            &         &    &    \\
       &    &    &     &            &   {\bf+1} &   2 &   1 &
&            &         &    &    \\
       &    &    &     &   {\bf--1} &   {\bf+4} &   6 &   4 &   1
&            &         &    &    \\
       &    &    &  1  &   {\bf--6} &  {\bf+15} &  20 &  15 &   6
&   {\bf--1} &         &    &    \\
       &    &  1 &  8  &  {\bf--28} &  {\bf+56} &  70 &  56 &  28
&   {\bf--8} &  {\bf+1}&    &    \\
       &  1 & 10 & 45  & {\bf--120} & {\bf+210} & 252 & 210 & 120
&  {\bf--45} & {\bf+10}&  1 &    \\ {\bf+1}& 12 & 66 & 220 &
{\bf--495} & {\bf+792} & 924 & 792 & 495 & {\bf--220} & {\bf+66}&
12 & 1

\end{tabular}

\centerline{$\cdots \cdots$}

\end{center}

Identity \eqref{e:122} corresponds to adding in each row the
elements marked with a plus sign and subtracting the ones marked
with a minus. By symmetry, we can represent this sum using only
the right-hand side of the above table, which by Theorem 2.1 is
the following Riordan array $\,\tilde{d}_{n,k}=\binom{2n}{n+k}$
with marked plus and minus entries

\begin{center}

\begin{tabular}{ccccccc}

  1 &         &          &            &          &    &     \\
  2 &  {\bf+1}&          &            &          &    &     \\
  6 &  {\bf+4}&  {\bf--1}&            &          &    &     \\
 20 & {\bf+15}&  {\bf--6}&   {\bf--1} &          &    &     \\
 70 & {\bf+56}& {\bf--28}&   {\bf--8} &  {\bf+1} &    &     \\
252 &{\bf+210}&{\bf--120}&  {\bf--45} & {\bf+10} &\ 1 &     \\
924 &{\bf+792}&{\bf--495}& {\bf--220} & {\bf+66} &\ 12&\ {\bf+1}

\end{tabular}

\centerline{$\cdots \cdots$}

\end{center}

\noindent In order to prove \eqref{e:122}, we wish to evaluate the
sum
\[
S_n=\sum_{j=-\infty}^\infty\bigg[\binom{2n}{n-5j-1}-\binom{2n}{n-5j-2}\bigg],
\]
i. e.,
\[
S_n=\sum_{k=0}^{\infty}\bigg[\binom{2n}{n\!+\!5k\!+\!1}-\binom{2n}{n\!+\!5k\!+\!2}-
\binom{2n}{n\!+\!5k\!+\!3}+\binom{2n}{n\!+\!5k\!+\!4}\bigg].
\]
In terms of the Riordan array
$\,\tilde{d}_{n,k}=\displaystyle\binom{2n}{n+k}$, $n,k\ge0$, we
have
\begin{equation}
S_n=\sum_{k=0}^\infty f_k\tilde{d}_{n,k}\,,
\label{e:s}
\end{equation}
where
\[
f(t)=\sum f_kt^k
=t-t^2-t^3+t^4+t^6-t^7-t^8+t^9+\cdots=\frac{t-t^2-t^3+t^4}{1-t^5}\,.
\]
For Pascal's triangle the generating function of the $A-$sequence
is $1+t$, since
$\binom{n+1}{k+1}=\binom{n}{k}+\binom{n}{k+1}\rule[-0.3cm]{0.cm}{0.6cm}$.
Hence, by Theorem 2.1, the generating function of the $A-$sequence
for $\{\tilde{d}_{n,k}\}$ is
\begin{equation}
A(t)=(1+t)^2.
\label{e:7}
\end{equation}
Either using Theorem 2.2 or noting that it follows from
\eqref{e:4} that $h(t)$ satisfies
\begin{equation}
t^2h^2+(2t-1)h+1=0,
\label{e:08}
\end{equation}
we obtain
$h(t)=\frac{1-2t-\sqrt{1-4t}}{2t^2}\rule[-0.25cm]{0.cm}{0.55cm}$.
Therefore, $\{\tilde{d}_{n,k}\}$ is the Riordan array characterized
by the pair $\big(d(t),h(t)\big)$, where
\[
d(t)=\sum_{n=0}^\infty\binom{2n}{n}t^n=\frac{1}{\sqrt{1-4t}}
\]
and
\begin{equation}
h(t)=\frac{1-2t-\sqrt{1-4t}}{2t^2}. \label{e:8}
\end{equation}
By Theorem 1.2, $S_n$ given by \eqref{e:s} satisfies
$S_n=[t^n]\,d(t)f(th(t))$. Note that
\[
f(t)=\frac{t\!-\!t^2\!-\!t^3\!+\!t^4}{1\!-\!t^5}=
\frac{t(1-t)(1-t^2)}{(1\!-\!t)(1\!+\!t\!+\!t^2\!+\!t^3\!+\!t^4)}=
\frac{t^{-1}-t}{t^{-2}+t^{-1}+1+t+t^2}\,.
\]
Setting $w:=th(t)$, by \eqref{e:08} $w$ and $w^{-1}$ are the roots
of the equation
\[
y^2+\frac{2t-1}{t}y+1=0.
\]
Hence,
\[
w^{-1}+w=\frac{1-2t}{t},\qquad
w^{-1}-w=\frac{\sqrt{1-4t}}{t}
\]
and
\[
w^2+w^{-2}=(w+w^{-1})^2-2=\frac{1-4t+2t^2}{t^2}.
\]
Therefore,
\[
f(th(t))=f(w)=\frac{\frac{\sqrt{1-4t}}{t}}
{\frac{1-4t+2t^2}{t^2}+\frac{1-2t}{t}+1}=
\frac{t\sqrt{1-4t}}{1-3t+t^2}
\]
and, finally,
\begin{equation}
\hspace{0.3cm} d(t)f\big(th(t)\big)=\frac{t}{1-3t+t^2}.
\rule[-0.5cm]{0.cm}{0.7cm}
\label{e:09}
\end{equation}
It is well known that both sequences of Fibonacci numbers $F_{2n}$
and $F_{2n+1}$ satisfy the recurrence relation
$x_n=3x_{n-1}-x_{n-2}$ and have generating functions
\[
\frac{t}{1-3t+t^2}=\sum_{n=0}^\infty F_{2n}t^n \hspace{0.5cm}
\]
and
\[
\frac{1-t}{1-3t+t^2}=\sum_{n=0}^\infty F_{2n+1}t^n,
\]
respectively (see \cite{14}, page 230). Hence, equation
\eqref{e:09} implies that
\[
d(t)f\big(th(t)\big)=\sum_{n=0}^\infty F_{2n}t^n, \hspace{0.37cm}
\]
from which \eqref{e:122} and \eqref{e:121} follow.

Identities \eqref{e:a3}, \eqref{e:a5}, and \eqref{e:a6} can be
obtained in a similar way. For \eqref{e:a3} and \eqref{e:a5} we
eliminate the even-numbered rows of Pascal's triangle

\begin{center}

\begin{tabular}{ccccc|ccccc}

       &       &       &       &     1 &   1 &    &    &   &    \\
       &       &       &     1 &     3 &   3 &  1 &    &   &    \\
       &       &     1 &     5 &    10 &  10 &  5 &  1 &   &    \\
       &$\cdot$&$\cdot$&$\cdot$&    35 &  35 & 21 &  7 & 1 &    \\
$\cdot$&$\cdot$&$\cdot$&$\cdot$&$\cdot$& 126 & 84 & 36 & 9 &  1

\end{tabular}

\centerline{$\cdots \cdots$}

\end{center}

\noindent and consider the right-hand side of what remains,
obtaining by Theorem 2.1 a Riordan array for which
\[
d(t)=\sum_{n=0}^\infty\binom{2n+1}{n+1}t^n=
\frac{1}{2t}\bigg(\frac{1}{\sqrt{1-4t}}-1\bigg)
\]
and
\[
h(t)=\frac{1-2t-\sqrt{1-4t}}{2t^2}.
\]

\vspace{0.6cm}

\section{Further Identities}

\vspace{0.2cm}

In this section, to illustrate the usefulness of the construction
considered in Theorem 2.1, we apply it to obtain a few more
identities via our Riordan array approach. Some of these
identities are well-known, while others are not. We believe that
identities \eqref{e:191} and \eqref{e:424} are new. We need one
more property of Riordan arrays, which generalizes a well-known
property of Pascal's triangle.

\vspace{0.5cm}

\noindent {\bf Theorem 4.1} {\em If $D=\big(d(t),h(t)\big)$ is a
Riordan array, then for any integers $k\ge s\ge1$ we have}
\begin{equation}
d_{n,k}=\sum_{j=s}^nd_{n-j,k-s}[t^{j}]\big(th(t)\big)^s.
\label{e:9}
\end{equation}
\label{t:7}

\noindent {\bf Proof.} From \eqref{e:1}, it follows that
$d_{n,k}=[t^n]\,d(t)\big(th(t)\big)^{k-s}\big(th(t)\big)^{s}$.
Hence,
\[
d_{n,k}=\sum_{j=s}^n\Big([t^{n-j}]\,d(t)\big(th(t)\big)^{k-s}\Big)
\Big([t^j]\big(th(t)\big)^s\Big),
\]
i.e.,
\[
d_{n,k}=\sum_{j=s}^nd_{n-j,k-s}[t^j]\big(th(t)\big)^s.\qquad\square
\]

\vspace{0.2cm}

One particular case of \eqref{e:9}, for $s=1$, is
\begin{equation}
d_{n,k}=\sum_{j=1}^nh_{j-1}d_{n-j,k-1},
\label{e:10}
\end{equation}
where $h(t)=\sum h_kt^k$. For example, in the case of Pascal's
triangle, since $h_k=1$ for all $k$, we obtain the well-known
identity
\begin{equation}
\sum_{j=1}^n\binom{n-j}{k-1}=\binom{n}{k}. \label{e:11}
\end{equation}
More generally, for Pascal's triangle, we have
\[
[t^j]\big(th(t)\big)^s=[t^j]\frac{t^s}{(1\!-\!t)^s}=
[t^{j-s}](1\!-\!t)^{-s}=\binom{-s}{j\!-\!s}(-1)^{j-1}=
\binom{j\!-\!1}{s\!-\!1}
\]
and, therefore, \eqref{e:9} becomes
\[
\sum_{j=s}^n\binom{n-j}{k-s}\binom{j-1}{s-1}= \binom{n}{k},
\]
which is (5.26) of \cite{19}, page 169, and contains \eqref{e:11}
as a particular case.

We now apply formulas \eqref{e:9} or \eqref{e:10} to Riordan
arrays obtained by the method described in Theorem 2.1.

\vspace{0.5cm}

\noindent {\bf Example 4.1} For fixed integers $p\ge2$ and
$r\ge0$, starting with Pascal's triangle and deleting $p-1$ rows
after each line kept, as described in Theorem 2.1, we obtain the
Riordan array
\[
d_{n,k}=\binom{pn+r}{(p-1)n+r+k}=\binom{pn+r}{n-k}.
\]
Note that
$d_{n,0}=\binom{pn+r}{(p-1)n+r}=\binom{pn+r}{n}\rule{0.cm}{0.5cm}$
and
\begin{equation}
d(t)=\sum_{n=0}^\infty\binom{pn+r}{n}t^n
=\sum_{n=0}^\infty\frac{(pn+r)!}{n!((p-1)n+r)!}t^n.
\label{e:23427}
\end{equation}
This is a hypergeometric series. By the same argument used for the
series \eqref{e:15h}, we get
\[
d(t)=\,_{p}F_{p-1}\left(
\begin{array}{c}
\frac{r+1}{p},\frac{r+2}{p},\ldots,\frac{r+p}{p}\\[4pt]
\frac{r+1}{p-1},\frac{r+2}{p-1},\ldots,\frac{r+p-1}{p-1}
\end{array}
\Bigg|\: \frac{p^pt}{(p\!-\!1)^{p-1}}\right).
\]
For $p=2$ and $r=0,1$, another expression in closed form for
\eqref{e:23427} is
\[
d(t)=\sum_{n=0}^\infty\binom{2n+r}{n}t^n
=\frac{1}{\sqrt{1-4t}}\bigg(\frac{1-\sqrt{1-4t}}{2t}\,\bigg)^{\!r}
\]
(see \cite{30}, Corollary 2.2).

In this example, the function $h$ is given by \eqref{e:23501} or,
equivalently, by \eqref{e:23421} with $q=p$. Combining \eqref{e:9}
and \eqref{e:1912} it follows that
\begin{equation}
\sum_{j=s}^n
\frac{ps}{(p-1)j+s}\binom{pj-1}{j-s}\binom{p(n-j)+r}{n-j-k+s}
=\binom{pn+r}{n-k}.
\label{e:237}
\end{equation}
In the particular case $s=1$, \eqref{e:237} becomes
\[
\sum_{j=1}^n
\frac{p}{(p-1)j+1}\binom{pj-1}{j-1}\binom{p(n-j)+r}{n-j-k+1}
=\binom{pn+r}{n-k},
\]
or,
\[
\sum_{j=1}^n\,
\frac{1}{pj+1}\binom{pj+1}{j}\binom{p(n-j)+r}{n-j-k+1}
=\binom{pn+r}{n-k}
\]
and, finally, adding $\binom{pn+r}{n-k+1}$ to both sides,
\begin{equation}
\sum_{j=0}^n\,
\frac{1}{pj+1}\binom{pj+1}{j}\binom{p(n-j)+r}{n-j-k+1}
=\binom{pn+r+1}{n-k+1}. \label{e:238}
\end{equation}
It would be nice to find a combinatorial interpretation for
\eqref{e:238}, since there are several interpretations for the
generalized Catalan number
$\frac{1}{pj+1}\binom{pj+1}{j}=\frac{1}{(p-1)j+1}\binom{pj}{j}$
(see \cite{19}, page 360, and \cite{29}).

Setting $j=i+s$, $x=ps$, $y=pk-ps+r$, and replacing $n$ by $n+k$,
identity \eqref{e:237} becomes
\[
\sum_{i=0}^n \frac{x}{x+pi}\binom{x+pi}{i}\binom{y+p(n-i)}{n-i}
=\binom{x+y+pn}{n},
\]
which is (5.62) of \cite{19}.

\vspace{0.5cm}

\noindent{\bf Example 4.2} We now consider Catalan's triangle
$d_{n,k}=\displaystyle\frac{k+1}{n+1}\binom{2(n+1)}{n-k}$, for
$n,k\ge0$,

\begin{center}

\begin{tabular}{cccccccc}

  1  &       &       &       &       &       &      &   \\
  2  &     1 &       &       &       &       &      &   \\
  5  &     4 &    1  &       &       &       &      &   \\
  14 &    14 &    6  &    1  &       &       &      &   \\
  42 &    48 &   27  &    8  &    1  &       &      &   \\
 132 &   165 &   110 &   44  &   10  &     1 &      &   \\
 429 &   572 &   429 &   208 &    65 &    12 &    1 &   \\

\end{tabular}
\centerline{$\cdots \cdots$}

\end{center}

This array was introduced in \cite{7} and has a nice
interpretation in terms of pairs of paths on a lattice. On the
bidimensional lattice $\mathbb{Z}^2$, consider all paths that
start at the origin, consist of unit steps and are such that all
steps go East or North. The length of a path is the number of
steps in the path. The distance between two paths of length $n$
with end-points $(a_n,b_n)$ and $(a_n',b_n')$, respectively, is
$|a_n-a_n'|$. Two paths are said to be non-intersecting if the
origin is the only point in common. Let $B(n,k)$, for $1\le k\le
n$, denote the number of pairs of non-intersecting paths of length
$n$ whose distance from one another is $k$. The array defined by
$d_{n,k}=B(n+1,k+1)$ is called Catalan's triangle. It is shown in
\cite{7} that $B(n,k)=\frac{k}{n}\binom{2n}{n-k}$ and that
\[
B(n-1,k-1)+2B(n-1,k)+B(n-1,k+1)=B(n,k).
\]
Therefore, by Theorem 1.3,
$d_{n,k}=\rule[-0.3cm]{0.cm}{0.5cm}\frac{k+1}{n+1}\binom{2(n+1)}{n-k}$
is a Riordan array and \eqref{e:7} is the generating function of
its $A-$sequence. The first column of this triangle is formed by
\[
d_{n,0}=\frac{1}{n+1}\binom{2(n+1)}{n}
=\frac{1}{n+2}\binom{2(n+1)}{n+1}=C_{n+1},
\]
where
$C_n=\frac{1}{n+1}\binom{2n}{n}=\binom{2n}{n}-\binom{2n}{n-1}$ is
the $n^{\rm th}$ Catalan number. We then have
\[
d(t)=h(t)=\frac{1-2t-\sqrt{1-4t}}{4t^2}.
\]
For fixed integers $p\ge2$ and $r\ge0$, we consider the Riordan
array
\begin{equation}
\tilde{d}_{n,k}=d_{pn+r,(p-1)n+k+r}=
\frac{(p-1)n+r+k+1}{pn+r+1}\binom{2(pn+r+1)}{n-k},
\label{e:1910}
\end{equation}
for which, by Theorem 2.1, the generating function of the
$A-$sequence is $A(t)=(1+t)^{2p}$. By Theorem 2.2, it follows that
for the Riordan array \eqref{e:1910} we have
\[
\big(h(t)\big)^s= \,_{2p}F_{2p-1}\left(
\begin{array}{c}
\frac{2ps}{2p},\frac{2ps\!+\!1}{2p},\ldots,\frac{2p(s\!+\!1)\!-\!1}{2p}\\[4pt]
\frac{2ps\!+\!1}{2p\!-\!1},\frac{2ps\!+\!2}{2p\!-\!1},\ldots,\frac{2p(s\!+\!1)\!-\!1}{2p\!-\!1}
\end{array}
\Bigg|\:\frac{(2p)^{2p}t}{(2p\!-\!1)^{2p-1}}\right)
\]
and
\begin{equation}
[t^j]\big(th(t)\big)^s=[t^{j-s}]\big(h(t)\big)^s=
\frac{2ps}{(2p-1)j+s}\binom{2pj-1}{j-s}.
\label{e:1914}
\end{equation}
From \eqref{e:9} and \eqref{e:1914} it follows that
\begin{equation}
\begin{split}
\sum_{j=s}^n\frac{2ps}{(2p\!-\!1)j\!+\!s}&\binom{2pj\!-\!1}{j\!-\!s}
\frac{(p\!-\!1)(n\!-\!j)\!+\!r\!+\!k\!-\!s\!+\!1}{p(n-j)+r+1}
\binom{2(p(n\!-\!j)\!+\!r\!+\!1)}{n-j-k+s}\\
&=\frac{(p\!-\!1)n+r+k+1}{pn+r+1}\binom{2(pn\!+\!r\!+\!1)}{n-k},
\rule{0.cm}{1.1cm}
\end{split}
\label{e:191}
\end{equation}
for every integers $s\le k\le n$, $p\ge1$, and $r\ge0$. Identity
\eqref{e:191} is probably new.

\vspace{0.5cm}

\noindent{\bf Example 4.3} In \cite{21} the following variant of
Catalan's triangle arises as an example of the infinite matrix
associated to a generating tree
\[
d_{n,k}=\left\{
\begin{array}{ll}
\displaystyle\frac{k+1}{n+1}\binom{2n-k}{n}, \ & \mbox{ if } \,
0\le k\le2n,  \\
0,   & \mbox{ if } \, 0\le2n<k. \end{array} \right.
\]

\begin{center}

\begin{tabular}{ccccccc}
  1 &     &     &     &     &     &    \\
  1 &   1 &     &     &     &     &    \\
  2 &   2 &   1 &     &     &     &    \\
  5 &   5 &   3 &  1  &     &     &    \\
 14 &  14 &   9 &  4  &  1  &     &    \\
 42 &  42 &  28 & 14  &  5  &  1  &    \\
132 & 132 &  90 & 48  & 20  &  6  & \ 1 \\
\end{tabular}
\centerline{$\cdots \cdots$}

\end{center}

\noindent In this case,
\[
\hspace{0.6cm}
d_{n+1,k+1}=d_{n,k}+d_{n,k+1}+d_{n,k+2}+\cdots+d_{n,n}
\]
and, thus,
\[
A(t)=1+t+t^2+\cdots=\frac{1}{1-t}.\hspace{0.5cm}
\]
From \eqref{e:4} it follows that $h(t)$ satisfies $\,th^2-h+1=0$
and, therefore, $h(t)=\frac{1-\sqrt{1-4t}}{2t}$. But this is
precisely the generating function of the Catalan numbers, which
form the first column of the triangle. Hence,
\[
d(t)=h(t)=\frac{1-\sqrt{1-4t}}{2t}.
\]
Note that
\[
\frac{1-\sqrt{1-4t}}{2t}=\,_2F_1\Bigg(
\begin{array}{c}
\frac12,\frac22\\[3pt]
\frac21
\end{array}\bigg|\:4t\Bigg).
\]
It follows from \eqref{e:235213} that
\[
\big(h(t)\big)^s=\,_2F_1\Bigg(
\begin{array}{c}
\frac{s}{2},\frac{s+1}{2}\\[2pt]
 s+1
\end{array}
\bigg|\:4t\Bigg)=
\sum_{n=0}^\infty\frac{s}{2n+s}\binom{2n+s}{n}t^n.
\]
For fixed integers $p\ge2$ and $r\ge0$, we consider the Riordan
array $\tilde{d}_{n,k}=d_{pn+r,(p-1)n+k+r}$. Then, for $n\ge0$,
\[
\tilde{d}_{n,k}=
\frac{(p-1)n+k+r+1}{pn+r+1}\binom{(p+1)n+r-k}{pn+r}
\]
if $\,0\le k\le(p+1)n+r$, and $\tilde{d}_{n,k}=0\,$ otherwise. By
Theorem 2.1, the generating function of the $A-$sequence of
$\{\tilde{d}_{n,k}\}$ is $\tilde{A}(t)=1/(1-t)^p$. By
\eqref{e:454545}, the $h-$function $\tilde{h}$ of
$\{\tilde{d}_{n,k}\}$ satisfies
\[
\begin{split}
[t^{n-1}]h(t)&=\frac{1}{n}[t^{n-1}]\big(A(t)\big)^n
=\frac{1}{n}[t^{n-1}](1-t)^{-pn}\\
&=\frac{(-1)^n}{n}\binom{-pn}{n-1}=\frac{(pn+n-2)!}{n!(pn-1)!}
\rule{0.cm}{1.cm}
\end{split}
\]
and, therefore,
\[
\tilde{h}(t)=\sum_{n=1}^\infty
\frac{1}{(p\!+\!1)n\!-\!1}\binom{(p\!+\!1)n\!-\!1}{n}t^{n-1}
=\sum_{n=0}^\infty
\frac{\big((p\!+\!1)n\!+\!p\!-\!1\big)!}{(n\!+\!1)!(pn\!+\!p\!-\!1)!}t^n.
\]
The function $\tilde{h}$ is hypergeometric. By the same argument
used for the series \eqref{e:15h}, we have
\[
\tilde{h}(t)=\,_{p+1}F_p\left(
\begin{array}{c}
\frac{p}{p+1},\frac{p+1}{p+1},\cdots,\frac{2p}{p+1}\\[4pt]
\frac{p+1}{p},\frac{p+2}{p},\cdots,\frac{2p}{p}
\end{array}
\Bigg|\:\frac{(p\!+\!1)^{p+1}}{p^p}t\right). \hspace{0.5cm}
\]
It follows from \eqref{e:235213} that
\[
\tilde{h}(t)=\left[ \,_{p+1}F_p\left(
\begin{array}{c}
\frac{1}{p+1},\frac{2}{p+1},\cdots,\frac{p+1}{p+1}\\[4pt]
\frac{2}{p},\frac{3}{p},\cdots,\frac{p+1}{p}
\end{array}
\Bigg|\:\frac{(p\!+\!1)^{p+1}}{p^p}t\right)\right]^p
\]
and also
\begin{equation}
\big(\tilde{h}(t)\big)^s=\,_{p+1}F_p\left(
\begin{array}{c}
\frac{sp}{p+1},\frac{sp+1}{p+1},\cdots,\frac{sp+p}{p+1}\\[4pt]
\frac{sp+1}{p},\frac{sp+2}{p},\cdots,\frac{sp+p}{p}
\end{array}
\Bigg|\: \frac{(p\!+\!1)^{p+1}}{p^p}t\right), \hspace{0.6cm}
\label{e:421}
\end{equation}
for any $s\in\mathbb{R}$. It is straightforward from \eqref{e:421}
that
\[
\big(\tilde{h}(t)\big)^s=\sum_{n=0}^\infty
\frac{sp}{(p+1)n+sp}\binom{(p+1)n+sp}{n}t^n
\]
and, therefore,
\[
[t^j]\big(t\tilde{h}(t)\big)^s=[t^{j-s}]\big(\tilde{h}(t)\big)^s=
\frac{sp}{(p+1)j-s}\binom{(p+1)j-s}{j-s}.
\]
Thus, by Theorem 4.1,
\[
\sum_{j=s}^{n-k+s}\!\!\!\frac{ps}{(p\!+\!1)j\!-\!s}\binom{\!(p\!+\!1)j\!-\!s}{j-s}\!
\frac{(p\!-\!1)(n\!-\!j)\!+\!k\!-\!s\!+\!r\!+\!1}{p(n\!-\!j)\!+\!r\!+\!1}
\binom{\!(p\!+\!1)(n\!-\!j)\!+\!r\!-\!k\!+\!s}{p(n\!-\!j)\!+\!r}
\]
\begin{equation}
=\frac{(p-1)n+k+r+1}{pn+r+1}\binom{(p+1)n+r-k}{pn+r}
\label{e:424}
\end{equation}
for $s\le k\le n$. Setting $j=i+s$, $x=ps$, $y=pk-ps+r$, and
replacing $n$ by $n+k$, identity \eqref{e:424} can be rewritten as
\begin{equation}
\begin{split}
\sum_{i=0}^{n}\frac{x}{(p\!+\!1)i\!+\!x}&\binom{(p\!+\!1)i\!+\!x}{i}
\frac{(p\!-\!1)(n\!-\!i)\!+\!y\!+\!1}{p(n\!-\!i)\!+\!y\!+\!1}
\binom{(p\!+\!1)(n\!-\!i)\!+\!y}{n-i}
\\
&\hspace{1.cm}=\frac{(p-1)n+x+y+1}{pn+x+y+1}\binom{(p+1)n+x+y}{n}
\rule{0.cm}{1.cm}.
 \end{split}
\label{e:450}
\end{equation}
Note that, for fixed $n\ge k\ge s\ge0$ and $r\ge0$ integers, our
argument only proves \eqref{e:450} for special values of $x$ and
$y$, namely of the form $x=ps$ and $y=p(k-s)+r$, with $p\ge2$ an
integer. As usual, this is enough to guarantee that \eqref{e:450}
holds for all $x$ and $y$ real, since both sides of \eqref{e:450}
are polynomials in $p$.

Identity \eqref{e:450} seems to be new, though it resembles the
old formula
\begin{equation}
\sum_{i=0}^{n}\frac{x}{zi\!+\!x}\binom{zi\!+\!x}{i}
\frac{y}{z(n\!-\!i)\!+\!y} \binom{z(n\!-\!i)\!+\!y}{n-i}=
\frac{x+y}{zn\!+\!x\!+\!y}\binom{zn\!+\!x\!+\!y}{n} \label{e:rh}
\end{equation}
due to Rothe (1793) and Hagen (1891) (see (3.142), in Gould's
collection \cite{31}, or (5.63) in, \cite{19}).

\vspace{0.6cm}

\section{Final comments}

\vspace{0.2cm}

Often identities of the same type as the ones obtained in Section
4 can be proved directly employing generating functions,
Lagrange's Inversion Formula, and standard Riordan array
techniques. Indeed, in \cite{32} Sprugnoli provides this kind of
proof to  most of the identities appearing in Gould's large
collection \cite{31}. As an example, we give a direct proof of
\eqref{e:450} along these lines. First note that
\[
\begin{split}
\frac{(p\!-\!1)n\!+\!y\!+1}{pn\!+\!y\!+\!1}\binom{(p\!+\!1)n\!+\!y}{n}&=
\bigg(1-\frac{n}{pn\!+\!y\!+\!1}\bigg)\binom{(p\!+\!1)n\!+\!y}{n}\\
=\binom{(p\!+\!1)n\!+\!y}{n}-\binom{(p\!+\!1)n\!+\!y}{n-1}
&=[t^n](1\!+\!t)^{(p+1)n+y}-[t^{n-1}](1\!+\!t)^{(p+1)n+y}\rule{0.cm}{0.95cm}\\
&=[t^n](1-t)(1+t)^y(1+t)^{(p+1)n}\rule{0.cm}{0.64cm}.
\end{split}
\]
On the other hand, by Lagrange's Inversion Formula \eqref{e:lag5},
\[
\begin{split}
\mathcal{G}\Big([t^n](1\!-\!t)(1\!+\!t)^y\big((1\!+\!t)^{p+1}\big)^n\Big)&=
\bigg[\frac{(1\!-\!w)(1\!+\!w)^y}{1-t(p\!+\!1)(1\!+\!w)^p}\,\bigg|\,
w=t(1\!+\!w)^{p+1}\bigg]\\
&=\bigg[\frac{(1\!-\!w)(1\!+\!w)^{y+1}}{1-pw}\,\bigg|\,w=t(1\!+\!w)^{p+1}\bigg].
\rule{0.cm}{0.9cm}
\end{split}
\]
Hence,
\begin{equation}
\begin{split}
\mathcal{A}(p,y;t)&:= \sum_{n=0}^\infty
\frac{(p\!-\!1)n\!+\!y\!+1}{pn\!+\!y\!+\!1}\binom{(p\!+\!1)n\!+\!y}{n}t^n\\
&=\bigg[\frac{(1\!-\!w)(1\!+\!w)^{y+1}}{1-pw}\bigg|w=t(1\!+\!w)^{p+1}\bigg].
\rule{0.cm}{0.9cm}
\end{split}
\label{e:460}
\end{equation}
By Theorem 2.1 of \cite{30}, we have the following identity for
the function given by \eqref{e:2352}
\begin{equation}
\begin{split}
\big(\mathcal{B}_{p+1}(t)\big)^x\!&=
\sum_{n=0}^\infty\frac{x}{(p\!+\!1)n+x}\binom{(p\!+\!1)n+x}{n}t^n\\
&=\big[(1\!+\!w)^x\big|w=t(1\!+\!w)^{p+1}\big].\rule{0.cm}{0.8cm}
\end{split}
\label{e:470}
\end{equation}
Combining \eqref{e:460} and \eqref{e:470} yields
\begin{equation}
\big(\mathcal{B}_{p+1}(t)\big)^x\cdot\mathcal{A}(p,y;t)
=\mathcal{A}(p,x+y;t)
\label{e:480}
\end{equation}
and, therefore, applying $[t^n]$ to both sides, \eqref{e:450}
holds. Identity of Rothe--Hagen \eqref{e:rh} mentioned above
follows from the application of $[t^n]$ to both sides of the
equality
$\big(\mathcal{B}_p(t)\big)^{x+y}=\big(\mathcal{B}_p(t)\big)^x\cdot\big(\mathcal{B}_p(t)\big)^y$.
It is interesting to observe that \eqref{e:480} implies that
$\mathcal{A}(p,y;t)$ is a function of exponential type on $y$, as
the ratio $\mathcal{A}(p,x+y;t)/\mathcal{A}(p,y;t)$ does not
depend on $y$. As a matter of fact, by \eqref{e:460} and
\eqref{e:470},
\[
\mathcal{A}(p,y;t)=\frac{(2-\mathcal{B}_{p+1}(t))\mathcal{B}_{p+1}(t)}
{1+p-p\mathcal{B}_{p+1}(t)}\big(\mathcal{B}_{p+1}(t)\big)^y.
\]
We can restate \eqref{e:480} in terms of hypergeometric functions.
If we consider the general term
\[
A_n=\frac{(p\!-\!1)n\!+\!y\!+1}{pn\!+\!y\!+\!1}
\binom{(p\!+\!1)n\!+\!y}{n}t^n
\]
of the power series in \eqref{e:460} and calculate the ratio of
consecutive terms, we find
\[
\frac{A_{n+1}}{A_n}=
\frac{\big(n\!+\!\frac{y+1}{p+1}\big)\big(n\!+\!\frac{y+2}{p+1}\big)\cdots\big(n\!+\!\frac{y+p+1}{p+1}\big)}
{\big(n\!+\!\frac{y+2}{p}\big)\big(n\!+\!\frac{y+3}{p}\big)\cdots\big(n\!+\!\frac{y+p+1}{p}\big)}\cdot
\frac{n\!+\!\frac{y+p}{p-1}}{n\!+\!\frac{y+1}{p-1}}\cdot
\frac{1}{n\!+\!1}\cdot\frac{(p+1)^{p+1}t}{p^p}\,,
\]
from which it follows that
\[
\mathcal{A}(p,y;t)=\,_{p+2}F_{p+1}\left(
\begin{array}{c}
\frac{y+1}{p+1},\frac{y+2}{p+1},\ldots,\frac{y+p+1}{p+1},\frac{y+p}{p-1}\\[4pt]
\frac{y+2}{p},\frac{y+3}{p},\ldots,\frac{y+p+1}{p},\frac{y+1}{p-1}
\end{array}
\Bigg| \: \frac{(p+1)^{p+1}t}{p^p}\right).
\]
Finally, identity \eqref{e:480} can be stated in terms of
hypergeometric functions as
\[
\,_{p+1}F_{p}\!\!\left(\!\!
\begin{array}{c}
\frac{x}{p+1}\!,\!\frac{x+1}{p+1}\!,\ldots,\!\frac{x+p}{p+1}\\[4pt]
\frac{x+1}{p}\!,\!\frac{x+2}{p}\!,\ldots,\!\frac{x+p}{p}
\end{array}
\!\Bigg|  \frac{(p\!+\!1\!)^{p+1}t}{p^p}\!\right)\!\cdot\!
\,_{p+2}F_{p+1}\!\!\left(\!\!
\begin{array}{c}
\frac{y+1}{p+1}\!,\!\frac{y+2}{p+1}\!,\ldots,\!\frac{y+p+1}{p+1}\!,\!\frac{y+p}{p-1}\\[4pt]
\frac{y+2}{p}\!,\!\frac{y+3}{p}\!,\ldots,\!\frac{y+p+1}{p}\!,\!\frac{y+1}{p-1}
\end{array}
\!\Bigg|  \frac{(p\!+\!1\!)^{p+1}t}{p^p}\!\right)
\]
\[
=\,_{p+2}F_{p+1}\!\left(\!\!\!
\begin{array}{c}
\frac{x+y+1}{p+1},\frac{x+y+2}{p+1},\ldots,\frac{x+y+p+1}{p+1},\frac{x+y+p}{p-1}\\[4pt]
\frac{x+y+2}{p},\frac{x+y+3}{p},\ldots,\frac{x+y+p+1}{p},\frac{x+y+1}{p-1}
\end{array}
\!\Bigg| \, \frac{(p\!+\!1)^{p+1}t}{p^p}\right).
\]
Of course this identity looks simpler if we further replace
$(p+1)^{p+1}t/p^p$ by $t$, but then the coefficients of the
corresponding developments in power series are no longer integers.

By a similar argument, \eqref{e:191} can be rephrased as
\begin{equation}
\begin{split}
\sum_{i=0}^n\frac{2x}{(2p\!-\!1)i\!+\!2x}&\binom{2pi\!+\!2x\!-\!1}{i}
\frac{(p\!-\!1)(n\!-\!i)\!+\!y\!+\!1}{p(n\!-\!i)\!+\!y\!+\!1}
\binom{2\big(p(n\!-\!i)\!+\!y\!+\!1\big)}{n-i}\\
&=\frac{(p\!-\!1)n\!+\!x\!+\!y\!+\!1}{pn\!+\!x\!+\!y\!+\!1}
\binom{2(pn\!+\!y\!+\!1)}{n}.\rule{0.cm}{0.9cm}
\end{split}
\label{e:520}
\end{equation}
Also as above, we find
\[
\begin{split}
\mathcal{C}(p,x;t)&:=
\sum_{n=0}^\infty\frac{2x}{(2p\!-\!1)n\!+\!2x}\binom{2pn\!+\!2x\!-\!1}{n}t^n\\
&=\big[(1+w)^{2x}\,\big|\,w=t(1+w)^{2p}\big]\rule{0.cm}{0.7cm}\\
\mathcal{D}(p,y;t)&:=
\sum_{n=0}^\infty\frac{(p\!-\!1)n\!+\!y\!+\!1}{pn\!+\!y\!+\!1}
\binom{2(pn\!+\!y\!+\!1)}{n}t^n\rule{0.cm}{0.9cm}\\
&=\bigg[\frac{(1-w)(1+w)^{2y+2}}{1+(1-2p)w}\bigg|w=t(1+w)^{2p}\bigg]
\rule{0.cm}{0.9cm}.
\end{split}
\]
Hence, $\mathcal{C}(p,x;t)=\big(\phi(t)\big)^x$ and
$\mathcal{D}(p,y;t)=\psi(t)\big(\phi(t)\big)^y$, where
$\phi(t)=\mathcal{C}(p,1;t)$ and
$\psi(t)=(2-\mathcal{C}(p,\frac{1}{2};t))\mathcal{C}(p,1;t)/
(2p+(1-2p)\mathcal{C}(p,\frac{1}{2};t))$. Therefore,
\begin{equation}
\mathcal{C}(p,x;t)\cdot\mathcal{D}(p,y;t)=\mathcal{D}(p,x+y;t).
\label{e:540}
\end{equation}
Applying $[t^n]$ to both sides, \eqref{e:520} follows. Using the
same argument as above, we can show that
\[
\mathcal{C}(p,x;t)=\!\,_{2p}F_{2p-1}\!\left(\!
\begin{array}{c}
\frac{2x}{2p},\frac{2x+1}{2p},\ldots,\frac{2x+2p-1}{2p}\\[4pt]
\frac{2x+1}{2p-1},\frac{2x+2}{2p-1},\ldots,\frac{2x+2p-1}{2p-1}
\end{array}
\!\Bigg| \, \frac{(2p)^{2p}\,t}{(2p\!-\!1)^{2p-1}}\right)=
\big(\mathcal{B}_{2p}(t)\big)^{2x}
\]
and
\[
\mathcal{D}(p,y;t)=\!\,_{2p+2}F_{2p+1}\!\left(\!
\begin{array}{c}
\frac{2y+3}{2p},\frac{2y+4}{2p},\ldots,\frac{2y+2p+2}{2p},\frac{y+p}{p-1},\frac{y+1}{p}\\[4pt]
\frac{2y+3}{2p-1},\frac{2y+4}{2p-1},\ldots,\frac{2y+2p+2}{2p-1},\frac{y+1}{p-1},\frac{y+p+1}{p}
\end{array}
\!\Bigg| \, \frac{(2p)^{2p}\,t}{(2p\!-\!1)^{2p-1}}\right).
\]
Hence, by \eqref{e:540}, another identity involving a product of
hypergeometric functions can be derived.

\vspace{0.5cm}

\noindent {\bf \large Acknowledgements.} The author wishes to
thank the referees for invaluable suggestions on the organization
of this paper, for providing reference \cite{30}, and for pointing
out that some identities obtained in this article (as
illustrations of our method) can alternately be obtained directly,
using the ideas in \cite{30}. We also express our gratitude to
Prof.\ J.\ Cigler for pointing out that a $q-$version of
\eqref{e:a1} and \eqref{e:a2} already appeared in the 1917 paper
\cite{27} by I.\ Schur.

\end{document}